\documentclass[11pt]{article}
\usepackage{mathrsfs}
\usepackage{amssymb}
\usepackage{amsmath}
\usepackage{amsbsy}
\usepackage{epsfig}
\usepackage{enumerate}
\usepackage{bm}

\usepackage{color}

\newcommand{\blue}[1]{\textcolor{blue}{#1}}

\usepackage[colorlinks,linkcolor=blue]{hyperref}

\newcommand{\ARXIV}[1]{\href{http://arXiv.org/abs/#1}{\blue{arXiv:#1}}}

\newcommand{\MR}[1]{\href{http://www.ams.org/mathscinet-getitem?mr=#1}{\blue{MR-#1}}}%

\def\proclaim#1{\par \smallskip\noindent {\bf #1}\bgroup\it\ }
\def\endproclaim{\egroup\par\smallskip}

\newtheorem{lemma}{Lemma}[section]
\newtheorem{theorem}{Theorem}[section]

\newtheorem{definition}{Definition}[section]
\newtheorem{remark}{Remark}[section]

\newbox\TempBox \newbox\TempBoxA

\def\Sbep{\widehat{\mathbb E}} 
\def\cSbep{\widehat{\mathcal E}} 

\def\Capc{\mathbb V} 
\def\cCapc{\mathcal V} 

\def\underwiggle 1{
\ifmmode\setbox\TempBox=\hbox{$ 1$}\else\setbox\TempBox=\hbox{
1}\fi \setbox\TempBoxA=\hbox to \wd\TempBox{\hss\char'176\hss}
\rlap{\copy\TempBox}\smash{\lower9pt\hbox{\copy\TempBoxA}} }
\renewcommand{\baselinestretch}{1.5}

\begin{document}

\thispagestyle{empty}

\begin{center}
 { \LARGE\bf Functional  central limit theorems  for random vectors under sub-linear expectations$^{\ast}$}
\end{center}

\begin{center} {\sc
\href{https://person.zju.edu.cn/en/stazlx}{\blue{Li-Xin Zhang}}\footnote{Research supported by grants from the NSF of China
(No.11731012),  Ten Thousands Talents Plan of Zhejiang Province (Grant No. 2018R52042)   and the Fundamental
Research Funds for the Central Universities..
}
}\\
{\sl \small School  of Mathematical Sciences, Zhejiang University, Hangzhou 310027} \\
(Email:stazlx@zju.edu.cn)    \\
\end{center}

\begin{abstract}The central limit theorem  of martingales is the fundamental tool for studying the convergence of stochastic processes. The central limit theorem and functional central limit theorem are obtained for martingale like random variables  under the sub-linear expectation by Zhang (2019).  In this paper, we consider the multi-dimensional martingale like random vectors and establish a functional central limit theorem. As applications, the Lindeberg central limit theorem for independent random vectors is established, and the sufficient and necessary conditions of the central limit theorem for independent and identically distributed random vectors are obtained.


\end{abstract}

\
\baselineskip 22pt

\renewcommand{\baselinestretch}{1.7}




\section{ Introduction and notations.}\label{sect1}
\setcounter{equation}{0}
Peng\cite{Peng07a}  introduced the notion of the sub-linear expectation. Under the sub-linear expectation, Peng\cite{Peng07a,Peng07b,Peng07c,Peng08a,Peng08b}  gave the notions of the G-normal
distributions, G-Brownian motions, G-martingales, independence of random variables, identical distribution of random variables and so on, and developed the weak law of large numbers and central limit theorem for independent and identically distributed (i.i.d.) random variables. Zhang\cite{Zh19} established the Lindeberg central limit theorem for independent but not necessary identically distributed one-dimensional random variables as well as martingale like sequences.  In this paper, we consider the multi-dimensional martingale like random vectors. In the classical probability space, since the convergence in distribution of a sequence of random vectors $\bm X_n=(X_{n,1},\ldots,X_{n,d})$ is equivalent to the convergence in distribution of any linear functions $\sum_k\alpha_kX_{n,k}$ of $\bm X_n$ by the Cram\'er-Wold device, the cental limit theorem for random vectors follows from   the cental limit theorem for one-dimensional random variables trivially. Under the sub-linear expectation, due to the non-linearity, the   Cram\'er-Wold device is no longer valid for showing the convergence of random vectors. In this paper, we  derive the functional Lindeberg central limit theorem for martingale like random vectors. As applications,
we  establish the Lindeberg central limit theorem for independent random vectors, give the sufficient and necessary conditions of the central limit theorem for independent and identically distributed random vectors,    obtain a L\'evy   characterization of a multi-dimensional G-Brownian motion and weaken a condition in the Lindeberg central limit theorems established by Zhang\cite{Zh19}.

We use the framework and notations of Peng \cite{Peng08b}. If the reader is familiar with these notations, the remainder of this section can be skipped.  Let  $(\Omega,\mathcal F)$
 be a given measurable space  and let $\mathscr{H}$ be a linear space of real functions
defined on $(\Omega,\mathcal F)$ such that if $X_1,\ldots, X_n \in \mathscr{H}$  then $\varphi(X_1,\ldots,X_n)\in \mathscr{H}$ for each
$\varphi\in C_{l,Lip}(\mathbb R^n)$,  where $C_{l,Lip}(\mathbb R^n)$ denotes the linear space of (local Lipschitz)
functions $\varphi$ satisfying
\begin{eqnarray*} & |\varphi(\bm x) - \varphi(\bm y)| \le  C(1 + |\bm x|^m + |\bm y|^m)|\bm x- \bm y|, \;\; \forall \bm x, \bm y \in \mathbb R^n,&\\
& \text {for some }  C > 0, m \in \mathbb  N \text{ depending on } \varphi. &
\end{eqnarray*}
$\mathscr{H}$ is considered as a space of ``random variables''. In this case, we denote $X\in \mathscr{H}$. We also denote   the   space of bounded Lipschitz
functions and the  space of bounded continuous functions on $\mathbb R^n$ by $C_{b,Lip}(\mathbb R^n)$ and  $C_b(\mathbb R^n)$, respectively.
 A  sub-linear expectation $\Sbep$ on $\mathscr{H}$  is a function $\Sbep: \mathscr{H}\to \overline{\mathbb R}$ satisfying the following properties: for all $X, Y \in \mathscr H$,
\begin{description}
  \item[\rm (1)]  Monotonicity: If $X \ge  Y$ then $\Sbep [X]\ge \Sbep [Y]$;
\item[\rm (2)] Constant preserving: $\Sbep [c] = c$;
\item[\rm (3)] Sub-additivity: $\Sbep[X+Y]\le \Sbep [X] +\Sbep [Y ]$ whenever $\Sbep [X] +\Sbep [Y ]$ is not of the form $+\infty-\infty$ or $-\infty+\infty$;
\item[\rm (4)] Positive homogeneity: $\Sbep [\lambda X] = \lambda \Sbep  [X]$, $\lambda\ge 0$.
 \end{description}
 Here $\overline{\mathbb R}=[-\infty, \infty]$. The triple $(\Omega, \mathscr{H}, \Sbep)$ is called a sub-linear expectation space. Give a sub-linear expectation $\Sbep $, let us denote the conjugate expectation $\cSbep$ of $\Sbep$ by
$  \cSbep[X]:=-\Sbep[-X]$,  $ \forall X\in \mathscr{H}$.
If $X$ is not in $\mathscr{H}$, we define its sub-linear expectation by $\Sbep^{\ast}[X]=\inf\{\Sbep[Y]:   X\le Y\in \mathscr{H}\}$.  When there is no ambiguity, we also denote it by  $\Sbep$.

After having the sub-linear expectation,
  we denote the pair $(\Capc,\cCapc)$ of capacities on  $(\Omega, \mathscr{H}, \Sbep)$ by setting
$$ \Capc(A):=\inf\{\Sbep[\xi]: I_A\le \xi, \xi\in\mathscr{H}\}, \;\; \cCapc(A):= 1-\Capc(A^c),\;\; \forall A\in \mathcal F, $$
where $A^c$  is the complement set of $A$.
It is obvious that $\Capc$ is sub-additive, i.e. $\Capc(A\bigcup B)\le \Capc(A)+\Capc(B)$. Also,
$\cCapc(A\bigcup B)\le \cCapc(A)+\Capc(B)$.

\bigskip
Next, we recall the notations of  identical distribution   and  independence.

\begin{definition} ({\em Pengc\cite{Peng07a,Peng08b}})

\begin{description}
  \item[ \rm (i)] ({\em Identical distribution}) Let $\bm X_1$ and $\bm X_2$ be two $n$-dimensional random vectors defined,
respectively, in sub-linear expectation spaces $(\Omega_1, \mathscr{H}_1, \Sbep_1)$
  and $(\Omega_2, \mathscr{H}_2, \Sbep_2)$. They are called identically distributed, denoted by $\bm X_1\overset{d}= \bm X_2$,  if
$$ \Sbep_1[\varphi(\bm X_1)]=\Sbep_2[\varphi(\bm X_2)], \;\; \forall \varphi\in C_{l,Lip}(\mathbb R^n), $$
whenever the sub-expectations are finite. A sequence $\{X_n;n\ge 1\}$ of random variables (or random vectors) is said to be identically distributed if $X_i\overset{d}= X_1$ for each $i\ge 1$.
\item[\rm (ii)] ({\em Independence})   In a sub-linear expectation space  $(\Omega, \mathscr{H}, \Sbep)$, a random vector $\bm Y =
(Y_1, \ldots, Y_n)$, $Y_i \in \mathscr{H}$ is said to be independent to another random vector $\bm X =
(X_1, \ldots, X_m)$ , $X_i \in \mathscr{H}$ under $\Sbep$,  if for each test function $\varphi\in C_{l,Lip}(\mathbb R^m \times \mathbb R^n)$
we have
$ \Sbep [\varphi(\bm X, \bm Y )] = \Sbep \big[\Sbep[\varphi(\bm x, \bm Y )]\big|_{\bm x=\bm X}\big],$
whenever $\overline{\varphi}(\bm x):=\Sbep\left[|\varphi(\bm x, \bm Y )|\right]<\infty$ for all $\bm x$ and
 $\Sbep\left[|\overline{\varphi}(\bm X)|\right]<\infty$.

 Random variables (or random vectors) $X_1,\ldots, X_n$ are said to be independent if for each $2\le k\le n$, $X_k$ is independent to $(X_1,\ldots, X_{k-1})$. A sequence of random variables (or random vectors) is said to be independent if for each $n$, $X_1,\ldots, X_n$ are independent.
\end{description}
\end{definition}

Finally, we recall the notations of G-normal distribution and  G-Brownian motion which are introduced by Peng\cite{Peng08b,Peng10}. We denote by $\mathbb S(d)$  the collection of all $d\times $ symmetric matrices. A function $G:\mathbb S(d) \to \mathbb R$ is called a sub-linear function monotonic in $A \in  \mathbb S(d)$ if  for each  $A, \overline{A}\in \mathbb S(d)$,
$$ \begin{cases} &G(A+\overline{A})\le G(A)+G(\overline{G}), \\
& G(\lambda A)=\lambda G(A), \;\; \forall \lambda>0,
\\
& G( A)\ge  G(\overline{A}), \;\; \text{ if } A\ge \overline{A}.
\end{cases}
$$
Here $A\ge \overline{A}$ means that $A- \overline{A}$ is semi-positive definitive.

\begin{definition} ({\em G-normal random variable}) Let $G:\mathbb S(d) \to \mathbb R$ be a continuous sub-linear function monotonic in $A \in  \mathbb S(d)$.
A $d$-dimensional random vector  $\bm \xi=(\xi_1,\ldots,\xi_d)$  in a sub-linear expectation space $(\widetilde{\Omega}, \widetilde{\mathscr H}, \widetilde{\mathbb E})$   is called a G-normal distributed  random variable   (written as $\xi\sim N\big(0, G\big)$  under $\widetilde{\mathbb E}$), if for any  $\varphi\in C_{l,Lip}(\mathbb R^d)$, the function $u(\bm x,t)=\widetilde{\mathbb E}\left[\varphi\left(\bm x+\sqrt{t} \bm \xi\right)\right]$ ($\bm x\in \mathbb R^d, t\ge 0$) is the unique viscosity solution of  the following heat equation:
      $$ \partial_t u -\frac{1}{2} G\left( D^2 u\right) =0, \;\; u(0,\bm x)=\varphi(\bm x). $$
where  $Du=\big(\partial_{x_i} u, i=1,\ldots,d\big)$ and $D^2u=D(Du)=\big(\partial_{x_i,x_j} u\big)_{i,j=1}^d$.
\end{definition}
That $\bm \xi$ is a G-normal distributed random vector  is equivalent to that,   if $\bm \xi^{\prime}$ is an independent copy of $\bm \xi$, then
$$ \widetilde{\mathbb E}\left[\varphi(\alpha \bm\xi+\beta \bm\xi^{\prime})\right]
=\widetilde{\mathbb E}\left[\varphi\big(\sqrt{\alpha^2+\beta^2}\bm \xi\big)\right], \;\;
\forall \varphi\in C_{l,Lip}(\mathbb R) \text{ and } \forall \alpha,\beta\ge 0, $$
and $G(A)=\widetilde{\mathbb E}\left[\langle\bm \xi A,\bm\xi\rangle\right]$ (cf. Definition II.1.4 and Example II.1.13 of Peng\cite{Peng10}), where $\langle\bm x,\bm y\rangle$ is the scalar product of $\bm x, \bm y$. When $d=1$, $G$ can be written as $G(\alpha)=\alpha^+\overline{\sigma}^2-\alpha^+\underline{\sigma}^2$, and we write $\xi\sim N(0,[\underline{\sigma}^2,\overline{\sigma}^2])$ if $\xi$ is a G-normal distributed random variable.

\begin{definition}\label{DefG-B} ({\em $G$-Brownian motion})  A $d$-dimensional random process $(\bm W_t)_{t\ge 0}$ in the sub-linear expectation space $(\widetilde{\Omega}, \widetilde{\mathscr H}, \widetilde{\mathbb E})$ is called a $G$-Brownian motion if
\begin{description}
  \item[\rm (i)] $\bm W_0=\bm 0$;
  \item[\rm (ii)]  For each $0\le t_1\le \ldots\le t_p\le t\le s$,
\begin{align}
&\widetilde{\mathbb E}\left[\varphi\big(\bm W_{t_1},\ldots, \bm W_{t_p}, \bm W_s-\bm W_t\big)\right]\nonumber
\\
= &
  \widetilde{\mathbb E}\left[\widetilde{\mathbb E}\left[\varphi\big(\bm x_1,\ldots, \bm x_p, \sqrt{t-s})\bm \xi\big)\right]\big|_{\bm x_1=\bm W_{t_1},\ldots, \bm x_p=\bm W_{t_p}}\right]
  \label{eqBrown} \\
  & \;\; \forall \varphi\in C_{l,Lip}(\mathbb R^{p\times(d+1)}), \nonumber
  \end{align}
  where $\bm \xi\sim N(0,G)$.
\end{description}
\end{definition}

Let $C_{[0,\infty)}=C_{[0,\infty)}(\mathbb R^d)$ be a function space of continuous real $d$-dimensional functions   on $[0,\infty)$ equipped with the supremum norm $\|\bm x\|=\sum\limits_{i=1}^{\infty}\sup\limits_{0\le t\le 2^i}|\bm x(t)|/2^i$, where $|\bm y|$ is the Euclidean norm of $\bm y$. Denote by $C_b\big(C_{[0,\infty)}\big)$  the set of bounded continuous  functions $h(x):C_{[0,\infty)}\to \mathbb R$. As showed in  Peng\cite{Peng08a,Peng10} and Denis, Hu, and Peng\cite{DenisHuPeng11},    there is a sub-linear expectation space $\big(\widetilde{\Omega}, \widetilde{\mathscr{H}},\widetilde{\mathbb E}\big)$ with
$\widetilde{\Omega}= C_{[0,\infty)}$ and $C_b\big(\widetilde{\Omega}\big)\subset \widetilde{\mathscr{H}}$ such that $(\widetilde{\mathscr{H}}, \widetilde{\mathbb E}[\|\cdot\|])$ is a Banach space, $\widetilde{\mathbb E}$ is countably sub-additive, and
the canonical process $W(t)(\omega) = \omega_t  (\omega\in \widetilde{\Omega})$ is a G-Brownian motion. In the sequel of this paper, the G-normal random vectors and G-Brownian motions are considered in $(\widetilde{\Omega}, \widetilde{\mathscr{H}}, \widetilde{\mathbb E})$.

\section{Functional Central limit theorem for martingale vectors.}\label{sectMain2}
\setcounter{equation}{0}

On the sub-linear expectation space $(\Omega, \mathscr{H}, \Sbep)$,  we write $\eta_n\overset{d}\to \eta$ if
$\Sbep\left[\varphi(\eta_n)\right]\to \Sbep\left[\varphi(\eta)\right]$ holds for all bounded and continuous functions $\varphi$, $\eta_n\overset{\Capc}\to \eta$ if $\Capc\left(|\eta_n-\eta|\ge \epsilon\right)\to 0$ for any $\epsilon>0$,   $  \eta_n\le  \eta +o(1)$ in capacity $\Capc$ if $ (\eta_n-\eta)^+\overset{\Capc}\to 0$,   $\eta_n\to   \eta $ in $L_p$ if $\lim_n \Sbep[|\eta_n-\eta|^p]=0$, and, $  \eta_n\le  \eta +o(1)$ in $L_p$ if $ (\eta_n-\eta)^+\to 0$ in $L_p$. We also write
$\xi\le \eta$ in $L_p$ if $\Sbep[((\xi-\eta)^+)^p]=0$, $\xi= \eta$ in $L_p$ if $\Sbep[|\xi-\eta|^p]=0$,  $X\le Y$ in   $\Capc$ if $\Capc\left(X-Y\ge \epsilon\right)=0$ for all $\epsilon>0$,
 and $X= Y$ in   $\Capc$ if   both $X\le Y$ and $Y\le X$ holds in $\Capc$.

 \begin{lemma}\label{lemma4.0.1} We have
 \begin{description}
   \item[\rm (1)]  if $X\le Y$ in $L_p$, then  $X\le Y$ in   $\Capc$;
   \item[\rm (2)]  if $X\le Y$ in   $\Capc$ and $\Sbep[((X-Y)^+)^p]<\infty$, then $X\le Y$ in $L_q$ for $0<q<p$;
   \item[\rm (3)]  if $X\le Y$ in   $\Capc$, $f(x)$ is non-decreasing continuous function  and $\Capc(|Y|\ge M)\to 0$ as $M\to \infty$, then $f(X)\le f(Y)$ in   $\Capc$;
   \item[\rm (4)] if $p\ge 1$, $X,Y\ge 0$ in $L_p$, $X\le Y$ in $L_p$, then $\Sbep[X^p]\le \Sbep[Y^p]$;
   \item[\rm (5)] if $\Sbep$ is countably additive, then $X\le Y$ in $\Capc$ is equivalent to  $X\le Y$ in $L_p$ for any $p>0$;
     \item[\rm (6)]  if $X_n\to  0$ in  $L_p$, then  $X_n\to  0$ in   $\Capc$ and in $L_q$ for $0<q<p$;
   \item[\rm (7)]  if $X_n\to  0$ in   $\Capc$ and $\Sbep[|X_n|^p]\le C<\infty$, then $X_n\to 0$ in $L_q$ for $0<q<p$.
 \end{description}
 \end{lemma}
Properties (1)-(5) are proved in Zhang\cite{Zh19}, and (6) and (7) can be proved in a similar way.

\bigskip

We  recall the definition of the conditional expectation under the sub-linear expectation. Let $(\Omega, \mathscr{H}, \Sbep)$ be a sub-linear expectation space.

  Let   $\mathscr{H}_{n,0}\subset \ldots\subset
 \mathscr{H}_{n,k_n}$ be subspaces of $\mathscr{H}$ such that
  \begin{description}
    \item[\rm (i)]  any constant   $c\in \mathscr{H}_{n,k}$ and,
    \item[\rm (ii)] if $X_1,\ldots,X_d\in \mathscr{H}_{n,k}$, then $\varphi(X_1,\ldots,X_d)\in \mathscr{H}_{n,k}$ for any $\varphi\in C_{l,lip}(\mathbb R^d)$, $k=0,\ldots, k_n$.
  \end{description}
  Denote $\mathscr{L}(\mathscr{H})=\{X:\Sbep[|X|]<\infty, X\in \mathscr{H}\}$.
 We consider a system of operators in $\mathscr{L}(\mathscr{H})$,
 $$\Sbep_{n,k}: \mathscr{L}(\mathscr{H})\to \mathscr{L}(\mathscr{H}_{n,k}) $$
 and denote $\Sbep[X|\mathscr{H}_{n,k}]=\Sbep_{n,k}[X]$, $\cSbep[X|\mathscr{H}_{n,k}]=-\Sbep_{n,k}[-X]$.  $\Sbep[X|\mathscr{H}_{n,k}]$ is called the conditional sub-linear expectation of $X$ given $\mathscr{H}_{n,k}$, $\Sbep_{n,k}$ is called the conditional expectation operator.
Suppose that the operators $\Sbep_{n,k}$ satisfy  the following properties:  for all $X, Y \in \mathscr{L}({\mathscr H})$,
\begin{description}
  \item[\rm (a)]   $ \Sbep_{n,k} [ X+Y]=X+\Sbep_{n,k}[Y]$ in $L_1$ if $X\in \mathscr{H}_{n,k}$, and $ \Sbep_{n,k} [ XY]=X^+\Sbep_{n,k}[Y]+X^-\Sbep_{n,k}[-Y]$ in $L_1$  if
$X\in \mathscr{H}_{n,k}$ and $XY\in \mathscr{L}({\mathscr H})$;
\item[\rm (b)]   $\Sbep\left[\Sbep_{n,k} [ X]\right]=\Sbep[X]$.
 \end{description}
In Zhang\cite{Zh19}, it has been shown  conditional expectation operators $\Sbep_{n,k}$ satisfy  that, for any $X,Y\in\mathscr{L}(\mathscr{H})$,
\begin{description}
  \item[\rm (c)]    $\Sbep_{n,k} [c] = c$ in $L_1$, $\Sbep_{n,k} [\lambda X] = \lambda \Sbep_{n,k}  [X]$ in $L_1$  if $\lambda\ge 0$;
  \item[\rm (d)]  $\Sbep_{n,k}[X]\le   \Sbep_{n,k}[Y]$ in $L_1$ if $X\le Y$ in $L_1$;
  \item[\rm (e)]  $\Sbep_{n,k}[X]-\Sbep_{n,k}[Y]\le \Sbep_{n,k}[X-Y]$ in $L_1$;
 \item[\rm (f)]  $\Sbep_{n,k}\left[\left[\Sbep_{n,l} [ X]\right]\right]=\Sbep_{n,l\wedge k} [ X]$ in $L_1$;
 \item[\rm (g)] if $|X|\le M$ in $L_p$ for all $p\ge  1$, then $ \big|\Sbep_{n,k}[X]\big| \le M$ in $L_p$ for all $p\ge  1$.
 \end{description}

For a random vector $\bm X=(X_1,\ldots, X_d)$, we denote $\Sbep[\bm X]=(\Sbep[X_1],\ldots, \Sbep[X_d])$ and $\Sbep[\bm X|\mathscr{H}_{n,k}]=(\Sbep[X_1|\mathscr{H}_{n,k}],\ldots, \Sbep[X_d|\mathscr{H}_{n,k}])$.
Now, we assume that $\{\bm Z_{n,k}; k=1,\ldots, k_n\}$ is an array of $d$-dimensional random vectors  such that $\bm Z_{n,k}\in \mathscr{H}_{n,k}$ and $\Sbep[|\bm Z_{n,k}|^2]<\infty$, $k=1,\ldots, k_n$.
Let $D_{[0,1]}=D_{[0,1]}(\mathbb R^d)$ be the space of right continuous $d$-dimensional functions having finite left limits which is endowed with the Skorohod topology (c.f. Billingsley\cite{B68}), $\tau_n(t)$ be a non-decreasing function in $D_{[0,1]}(\mathbb R^1)$ which takes integer values with $\tau_n(0)=0$, $\tau_n(1)=k_n$. Define  $\bm S_{n,i}=\sum_{k=1}^i \bm Z_{n,k}$,
\begin{equation}\label{eqthFCLTM.1} \bm W_n(t)= \bm S_{n, \tau_n(t)}.
\end{equation}
Then $\bm W$ is element in $D_{[0,1]}(\mathbb R^d)$.  The following is the functional central limit theorem.

 \begin{theorem} \label{thFCLTM}   Suppose that the operators $\Sbep_{n,k}$ satisfy (a) and (b).   Assume that   the following Lindeberg condition is satisfied:
 \begin{equation}\label{eqLindebergM}
  \sum_{k=1}^{k_n}\Sbep\left[ \left( |\bm Z_{n,k}|^2-\epsilon \right)^+|\mathscr{H}_{n,k-1}\right]\overset{\Capc}\to 0\;\; \forall \epsilon>0,
 \end{equation}
 and
  \begin{equation}\label{eqCLTCondM.2}
  \sum_{k=1}^{k_n}\left\{ |\Sbep[\bm Z_{n,k} |\mathscr{H}_{n,k-1}]|+|\cSbep[\bm Z_{n,k} |\mathscr{H}_{n,k-1}]|\right\}    \overset{\Capc}\to 0.
 \end{equation}
  Further,   assume that there is  a continuous non-decreasing non-random function   $\rho(t)$ and a non-random function $G:\mathbb S(d)\to \mathbb R$ for which
 \begin{equation}\label{eqFCLTCondM.3}
 \sum_{k\le \tau_n(t)}  \Sbep\left[\langle \bm Z_{n,k} A,\bm Z_{n,k} \rangle \big|\mathscr{H}_{n,k-1}\right]    \overset{\Capc}\to G(A)\rho(t), \;\; A\in \mathbb S(d).
 \end{equation}
 Then  for any $0=t_0<\ldots< t_d\le 1$,
 \begin{equation} \label{eqfinitedimension} \Big( \bm W_n(t_1),\ldots, \bm W_n(t_d)\Big)\overset{d}\to \Big( \bm W(\rho(t_1)),\ldots, \bm W(\rho(t_d))\Big),
 \end{equation}
 and for any  bounded continuous function $\varphi:D_{[0,1]}(\mathbb R^d)\to \mathbb R$,
  \begin{equation} \label{eqFCLTM} \lim_{n\to \infty}\Sbep\left[\varphi\left(\bm W_n\right)\right]=\widetilde{\mathbb E}[\varphi(\bm W\circ\rho)],
  \end{equation}
   where  $\bm W$ is $G$-Brownian motion with  $\bm W(1) \sim N(0,G)$  under $\widetilde{\mathbb E}$, and $\bm W\circ\rho(t)=\bm W(\rho(t))$.
 \end{theorem}

\begin{remark}\label{remark2.1}
Let $G_n(A,t)=\sum_{k\le \tau_n(t)}  \Sbep\left[\langle \bm Z_{n,k} A,\bm Z_{n,k} \rangle \big|\mathscr{H}_{n,k-1}\right]$. It is easily seen that
$G_n(A,t):\mathbb S(d) \to \mathbb R$ be a continuous sub-linear function monotonic in $A \in  \mathbb S(d)$. So, $G$ is a continuous sub-linear function monotonic in $A \in  \mathbb S(d)$. Without loss of generality, we assume $G(I_{d\times d})=1$ for otherwise we can replace $\rho(t)$ by $G(I_{d\times d})\rho(t)$. It is obvious that
\begin{align*}
|G_n(A,t)-G_n(\overline{A},t)|\le & d\|A-\overline{A}\|_{\infty}
\sum_{k\le \tau_n(t)}  \Sbep\left[\langle \bm Z_{n,k},\bm Z_{n,k} \rangle \big|\mathscr{H}_{n,k-1}\right]\\
 = & d\|A-\overline{A}\|_{\infty} G_n(I,t).
 \end{align*}
It follows that $|G(A)-G(\overline{A})|\le   d\|A-\overline{A}\|_{\infty}$. Then, it can be verified that (\ref{eqFCLTCondM.3}) holds uniformly in $A$ in a bounded area, and $G(A)$ is continuous in $A\in \mathbb S(d)$.
\end{remark}

The proof of this theorem will stated in the last section.
\begin{remark}
When $d=1$, (\ref{eqFCLTCondM.3}) is equivalent to
\begin{equation}\label{eqremark2.1}
 \sum_{k\le \tau_n(t)}  \Sbep[Z_{n,k}^2|\mathscr{H}_{n,k-1}]    \overset{\Capc}\to \rho(t), \;\; t\in[0,1],
 \end{equation}
 \begin{equation}\label{eqremark2.2}
\text{and}\; \sum_{k\le \tau_n(t)}  \cSbep[Z_{n,k}^2|\mathscr{H}_{n,k-1}]    \overset{\Capc}\to r\rho(t), \;\; t\in[0,1].
 \end{equation}
The condition (\ref{eqremark2.1}) is assumed in Zhang\cite{Zh19}. But, (\ref{eqremark2.2}) is replaced by a more stringent  condition as follows,
$$
 \sum_{k=1}^{k_n} \left| r\Sbep[Z_{n,k}^2|\mathscr{H}_{n,k-1}] - \cSbep[Z_{n,k}^2|\mathscr{H}_{n,k-1}]\right|  \overset{\Capc}\to 0.
$$
As shown in Remark \ref{remark3.1}, (\ref{eqremark2.1}) and (\ref{eqremark2.2}) can not be weakened furthermore.
\end{remark}


\section{Applications}\label{Appl}
\setcounter{equation}{0}
From Theorem \ref{thFCLTM}, we have the following functional central limit theorem for independent random vectors.

\begin{theorem} \label{th2}   Let  $\{\bm X_{n,k};k=1,\ldots, k_n\}$   be an array of independent $d$-dimensional random vectors,  $n=1,2,\ldots$, $\tau_n(t)$ be a non-decreasing function in $D_{[0,1]}(\mathbb R^1)$ which takes integer values with $\tau_n(0)=0$, $\tau_n(1)=k_n$.
Denote
$ \bm W_n(t)=\sum_{k\le \tau_n(t)}\bm X_{n,k}. $
 Assume that
 \begin{equation} \label{eqth2.2}
  \sum_{k=1}^{k_n}\Sbep\left[ \left( |\bm X_{n,k}|^2-\epsilon \right)^+ \right] \to 0\;\; \forall \epsilon>0,
 \end{equation}
 and
  \begin{equation} \label{eqth2.3}
  \sum_{k=1}^{k_n}\left\{ |\Sbep[\bm X_{n,k}  ]|+|\cSbep[\bm X_{n,k}  |\right\}     \to 0.
 \end{equation}
Further,   assume that there is  a continuous non-decreasing non-random function   $\rho(t)$ and a non-random function $G:\mathbb S(d)\to \mathbb R$ for which
 \begin{equation}\label{eqth2.4}
 \sum_{k\le \tau_n(t)}  \Sbep\left[\langle \bm X_{n,k} A,\bm X_{n,k} \rangle\right]    \overset{\Capc}\to G(A)\rho(t), \;\; A\in \mathbb S(d).
 \end{equation}
 Then  for any $0=t_0<\ldots< t_d\le 1$,
 \begin{equation}  \label{eqth2.6} \Big( \bm W_n(t_1),\ldots, \bm W_n(t_d)\Big)\overset{d}\to \Big( \bm W(\rho(t_1)),\ldots, \bm W(\rho(t_d))\Big),
 \end{equation}
 and for  any   continuous function $\varphi:D_{[0,1]}(\mathbb R^d)\to \mathbb R$ with $ |\varphi(\bm x)|\le C \sup_{t\in [0,1]}|\bm x(t)|^2$,
  \begin{equation} \label{eqth2.7}  \lim_{n\to \infty}\Sbep\left[\varphi\left(\bm W_n\right)\right]=\widetilde{\mathbb E}[\varphi(\bm W\circ\rho)],
  \end{equation}
   where  $\bm W$ is $G$-Brownian motion on $[0,\infty)$ with  $\bm W(1) \sim N(0,G)$  under $\widetilde{\mathbb E}$. Further,  when $p>2$, (\ref{eqth2.7}) holds  for any   continuous function $\varphi:D_{[0,1]}(\mathbb R^d)\to \mathbb R$ with $ |\varphi(\bm x)|\le C \sup_{t\in [0,1]}|\bm x(t)|^p$ if    (\ref{eqLindeberg}) is replaced by the condition that
 \begin{align}\label{eqth2.8}
 \sum_{k=1}^{k_n}\Sbep\left[ \left| \bm X_{n,k}  \right|^p\right]\to 0.
 \end{align}
\end{theorem}

{\bf Proof.}
For a bounded continuous function $\varphi$, (\ref{eqth2.7}) follows from  Theorem \ref{thFCLTM}  for the functional central limit theorem of martingale vectors. For continuous function $\varphi:D_{[0,1]}(\mathbb R^d)\to \mathbb R$ with $ |\varphi(\bm x)\le C \sup_{t\in [0,1]}|\bm x(t)|^p$, we first note that (\ref{eqth2.2}) is implied by  (\ref{eqth2.8}) for $p>2$. Since  (\ref{eqth2.7}) holds for bounded continuous function $\varphi$ and
$$\big|\varphi(\bm x)-(-N)\vee \varphi(x)\wedge N\big|\le \big(C\sup_{t\in [0,1]}|\bm x(t)|^p-N)^+,  $$
 it is sufficient to show that $\{\max\limits_{i\le k_n}|\sum\limits_{k\le i}\bm X_{n,k}|^p, n\ge 1\}$ is uniformly integrable,
 i.e.,
\begin{equation}\label{eqth2proof.2}
\lim_{N\to \infty}\limsup_{n\to \infty}\Sbep\left[ \left(  \max_{i\le k_n}\Big|\sum_{k=1}^i \bm X_{n,k}\Big|^p -N\right)^+\right]=0
\end{equation}
under the conditions (\ref{eqth2.2}), (\ref{eqth2.4}),  (\ref{eqth2.3}) or/and (\ref{eqth2.8}).  For showing (\ref{eqth2proof.2}), it is sufficient to consider the one-dimensional case.  Let $Y_{n,k}=(-  1)\vee X_{n,k}\wedge 1$ and   $\widehat{Y}_{n,k}=X_{n,k}-Y_{n,k}$. Then,  the Lindeberg condition (\ref{eqth2.2})  implies that
\begin{equation}\label{eqth2proof.3}
   \sum_{k=1}^{k_n} \Sbep[|\widehat{Y}_{n,k} |]
=   \sum_{k=1}^{k_n}\Sbep\Big[\big(|X_{n,k}|- 1\big)^+\Big]
\le  2 \sum_{k=1}^{k_n}\Sbep\Big[\big(|X_{n,k}|^2-  1/2\big)^+\Big]  \to 0.
\end{equation}
 It follows that
  \begin{equation}\label{eqth2proof.4}
 \sum_{k=1}^{k_n}\left\{ |\Sbep[Y_{n,k}]|+|\cSbep[Y_{n,k}]|\right\}   \to 0,
 \end{equation}
by (\ref{eqth2.3}).  Also, it is obvious that
  \begin{equation}\label{eqth2proof.5} \sum_{k=1}^{k_n} \Sbep[|Y_{n,k}|^q]   \le \sum_{k=1}^{k_n} \Sbep[Y_{n,k}^2]\le \sum_{k=1}^{k_n} \Sbep[X_{n,k}^2]  =O(1), \forall q\ge 2.
   \end{equation}
  By the Rosenthal-type inequality for independent random variables (c.f. Theorem 2.1 of Zhang\cite{Zh16}),
   \begin{align}\label{eqth2proof.6}
   \Sbep\left[ \max_{i\le k_n}\Big|\sum_{k=1}^i Y_{n,k}\Big|^q\right]
   \le & C_q\left\{\sum_{k=1}^{k_n}\Sbep\left[|Y_{n,k}|^q\right]+\left(\sum_{k=1}^{k_n}\Sbep\left[Y_{n,k}^2\right]\right)^{q/2}\right. \nonumber \\
  &\quad + \left. \left(\sum_{k=1}^{k_n}\Big(\big|\Sbep\left[ Y_{n,k} \right]\big|+\big|\cSbep\left[ Y_{n,k} \right]\big|\Big)\right)^q\right\}
  \le    C_q,
   \end{align}
by (\ref{eqth2proof.4}) and (\ref{eqth2proof.5}). It follows that
\begin{align*}
&\lim_{N\to \infty}\limsup_{n\to \infty}\Sbep\left[ \left(  \max_{i\le k_n}\Big| \sum_{k=1}^i Y_{n,k} \Big|^p -N\right)^+\right]\\
&\; \le
\lim_{N\to \infty}\limsup_{n\to \infty}N^{-2}\Sbep\left[   \max_{i\le k_n}\Big| \sum_{k=1}^i Y_{n,k} \Big|^{2p} \right]= 0.
\end{align*}
For $\widehat{Y}_{n,k}$, by  the Rosenthal-type inequality for independent random variables again we have
\begin{align*}
&\Sbep\left[ \max_{i\le k_n}\Big|\sum_{k=1}^i \widehat{Y}_{n,k}\Big|^{p}\right] \\
\le &C_p\left\{  \sum_{k=1}^{k_n}\Sbep[|\widehat{Y}_{n,k}|^p] +\left( \sum_{k=1}^{k_n}\Sbep[|\widehat{Y}_{n,k}|^2] \right)^{p/2}
 + \left( \sum_{k=1}^{k_n}\big((\Sbep[\widehat{Y}_{n,k}])^+ +(\cSbep[\widehat{Y}_{n,k}])^-\big)\right)^{p }\right\}\\
 \le &C_p\left\{  \sum_{k=1}^{k_n}\Sbep[(|X_{n,k}|^p]-  1)^+ +\left( \sum_{k=1}^{k_n}\Sbep[(X_{n,k}^2-1)^+] \right)^{p/2}
  + \left( \sum_{k=1}^{k_n}\Sbep[(|X_{n,k}|-  1)^+] \right)^{p }\right\}\\
 &\to 0
\end{align*}
by  (\ref{eqth2proof.3}) and the condition (\ref{eqth2.2}) (and  (\ref{eqth2.8}) when $p>2$). Hence, (\ref{eqth2proof.2}) is proved. \hfill $\Box$

\begin{remark}\label{remark3.1} When $d=1$, the condition (\ref{eqth2.4}) is equivalent to
\begin{equation}\label{eqrk3.1.1}
 \sum_{k\le \tau_n(t)}  \Sbep[X_{n,k}^2 ]    \to \rho(t), \;\; t\in[0,1],
\end{equation}
\begin{equation}\label{eqrk3.1.2}
 \sum_{k\le \tau_n(t)}  \cSbep[X_{n,k}^2 ]     \to r\rho(t), \;\; t\in[0,1].
\end{equation}
 Suppose that $\{X_{n,k};k=1,\ldots, k_n\}$  is an array of independent random variables with $\Sbep[X_{n,k}]=\cSbep[X_{n,k}]=0$, $k=1,\ldots, k_n$, and the Lindeberg condition (\ref{eqLindeberg0}) is satisfied. If (\ref{eqth2.6}) or (\ref{eqth2.7}) holds, then as   shown in the proof of Theorem \ref{th2},
$$ \sum_{k\le k_n(t)}\Sbep[X_{n,k}^2]=\Sbep[W_n^2(t)] \to \Sbep[W^2(\rho(t))]=\rho(t), $$
$$\sum_{k\le k_n(t)}\cSbep[X_{n,k}^2]=\cSbep[W_n^2(t)] \to \cSbep[W^2(\rho(t))]=r\rho(t). $$
So, the conditions (\ref{eqrk3.1.1}) and (\ref{eqrk3.1.2})  can not be weakened furthermore.
\end{remark}

 Zhang\cite{Zh19}  gave the following Lindeberg's  central limit theorem for arrays of independent random variables.

\begin{proclaim}{Theorem A}  Let  $\{X_{n,k};k=1,\ldots, k_n\}$   be an array of independent random variables,  $n=1,2,\ldots$.
Denote    $\overline{\sigma}_{n,k}^2=\Sbep[X_{n,k}^2]$,  $\underline{\sigma}_{n,k}^2=\cSbep[X_{n,k}^2]$ and $B_n^2=\sum_{k=1}^{k_n} \overline{\sigma}_{n,k}^2$. Suppose that the Lindeberg condition is satisfied:
 \begin{align}\label{eqLindeberg0}
 \frac{1}{B_n^2}\sum_{k=1}^{k_n}\Sbep\left[ \left( X_{n,k}^2-\epsilon B_n^2 \right)^+\right]\to 0\;\; \forall \epsilon>0,
 \end{align}
 and further, there is a constant $r\in[0,1]$ such that
  \begin{align}\label{eqCLTCond0.2}
 \frac{\sum_{k=1}^{k_n} \left| r\overline{\sigma}_{n,k}^2 - \underline{\sigma}_{n,k}^2\right| }{B_n^2}\to 0, \;\; \text{also, }
 \end{align}
  \begin{align}\label{eqCLTCond0.3}
 \frac{\sum_{k=1}^{k_n}\left\{ |\Sbep[X_{n,k}]|+|\cSbep[X_{n,k}]|\right\}  }{B_n}\to 0.
 \end{align}
 Then for any  bounded continuous function $\varphi$,
  \begin{equation} \label{eqCLT0} \lim_{n\to \infty}\Sbep\left[\varphi\left(\frac{\sum_{k=1}^{k_n}X_{n,k}}{B_n}\right)\right]=\widetilde{\mathbb E}[\varphi(\xi)],
  \end{equation}
   where   $\xi\sim N(0,[r, 1])$  under $\widetilde{\mathbb E}$.
 \end{proclaim}

 Zhang\cite{Zh19}  also showed that the condition (\ref{eqCLTCond0.2}) can not be weakened to
\begin{equation}\label{eqCLTCond0.4}  \frac{\;\;\sum_{k=1}^{k_n}  \underline{\sigma}_{n,k}^2\;\;}{\sum_{k=1}^{k_n}  \overline{\sigma}_{n,k}^2} \to r.
\end{equation}
The following theorem  shows that if  we consider  a   sequence   of  independent random variables instead of arrays of independent random variables, then the condition (\ref{eqCLTCond0.2}) can   be weakened to  (\ref{eqCLTCond0.4}).
\begin{theorem}\label{thCLT}  Let   $\{X_k;k=1,2,\ldots\}$ be a sequence of independent random variables.
Denote    $\overline{\sigma}_{k}^2=\Sbep[X_{k}^2]$,  $\underline{\sigma}_{k}^2=\cSbep[X_{k}^2]$, $B_n^2=\sum_{k=1}^{n} \overline{\sigma}_{k}^2$ .  Suppose that the Lindeberg condition is satisfied:
 \begin{align}\label{eqLindeberg}
 \frac{1}{B_n^2}\sum_{k=1}^n\Sbep\left[ \left( X_{k}^2-\epsilon B_n^2 \right)^+\right]\to 0\;\; \forall \epsilon>0,
 \end{align}
 and further, there is a constant $r\in[0,1]$ such that
  \begin{align}\label{eqCLTCond.2}
   \frac{\;\;\sum_{k=1}^n  \underline{\sigma}_k^2\;\;}{\sum_{k=1}^{n}  \overline{\sigma}_k^2} \to r, \;\; \text{also, }
 \end{align}
  \begin{align}\label{eqCLTCond.3}
 \frac{\sum_{k=1}^{n}\left\{ |\Sbep[X_k]|+|\cSbep[X_k]|\right\}  }{B_n}\to 0.
 \end{align}
 Then for any   continuous function $\varphi$ with $\ |\varphi(x)|\le C x^2$,
  \begin{equation} \label{eqCLT} \lim_{n\to \infty}\Sbep\left[\varphi\left(\frac{\sum_{k=1}^n X_k}{B_n}\right)\right]=\widetilde{\mathbb E}[\varphi(\xi)],
  \end{equation}
   where   $\xi\sim N(0,[r, 1])$  under $\widetilde{\mathbb E}$. Further, when $p>2$, (\ref{eqCLT}) holds for any   continuous function   $\varphi$ with $\ |\varphi(x)|\le C |x|^p$    if (\ref{eqLindeberg}) is replaced by the condition that
 \begin{align}\label{eqLindeberg2}
 \frac{1}{B_n^p}\sum_{k=1}^n\Sbep\left[ \left| X_k  \right|^p\right]\to 0.
 \end{align}
 \end{theorem}

{\bf Proof}.
For proving Theorem \ref{thCLT},  we let  $k_n=n$, $X_{n,k}=X_k/B_n$, $k=1,\ldots, n$. It is easily seen that the array $\{X_{n,k}; k=1,\ldots, k_n\}$ satisfies (\ref{eqth2.2}) and (\ref{eqth2.3}). Denote $B_0=0$. Define the function $\tau_n(t)$ by
$$ \tau_n(t)=k \text{ if } B_k^2/B_n^2\le t<B_{k+1}^2/B_n^2, \;\; \text{ and } \tau_n(1)=n. $$
From the the Lindeberg condition (\ref{eqLindeberg0}), it is easily verified that
$$ \frac{\max_k \underline{\sigma}_k^2}{B_n^2}\le \frac{\max_k \overline{\sigma}_k^2}{B_n^2}\to 0. $$
It follows that
$$ \Big|\sum_{k\le \tau_n(t)}\Sbep[X_{n,k}^2]-t\Big|=\Big|\frac{B^2_{\tau_n(t)}}{B_n^2}- t\Big|\le \frac{\max_k \overline{\sigma}_k^2}{B_n^2}\to 0,  $$
and $\tau_n(t)\to \infty$ if $t>0$. By the condition (\ref{eqCLTCond.2}), we have
$$ \sum_{k\le \tau_n(t)}\cSbep[X_{n,k}^2]=\frac{\sum_{k\le \tau_n(t)} \underline{\sigma}_k^2}{B_n^2} = \frac{\sum_{k\le \tau_n(t)} \underline{\sigma}_k^2}{\sum_{k\le \tau_n(t)} \overline{\sigma}_k^2}\frac{B^2_{\tau_n(t)}}{B_n^2} \to rt.  $$
So, (\ref{eqrk3.1.1}) and (\ref{eqrk3.1.2}) are satisfied with $\rho(t)=t$.  Hence, (\ref{eqCLT}) follows from   (\ref{eqth2.7}). \hfill $\Box$

The next theorem gives the sufficient and necessary conditions of the central limit theorem for independent and identically distributed random vectors.
For a one-dimensional random variable $X$ and a positive constant $c$, we write $X^{(c)}=(-c)\vee(X\wedge c)$, and for a random vector $\bm X=(X_1,\ldots, X_d)$ we write
$\bm X^{(c)}=(X_1^{(c)},\ldots, X_d^{(c)})$.

\begin{theorem} \label{cltiid} Let   $\{\bm X_k;k=1,2,\ldots\}$ be a sequence of independent and identically distributed $d$-dimensional  random vectors, $\bm S_n=\sum_{k=1}^n \bm X_k$. Suppose that
  \begin{description}
    \item[\rm (i) ] $\lim\limits_{c\to\infty} \Sbep[|\bm X_1|^2\wedge c]$ is finite;
    \item[\rm (ii)]  $x^2\Capc\left(|\bm X_1|\ge x\right)\to 0$ as $x\to \infty$;
    \item[\rm (iii)]  $\lim\limits_{c\to \infty}\Sbep\left[  \bm X_1^{(c)}\right]=\lim\limits_{c\to \infty}\Sbep\left[-\bm X_1^{(c)}\right]=\bm 0$;
     \item[\rm (iv)] The limit
     \begin{equation} \label{cltiiid1} G(A)=\lim_{c\to\infty}\Sbep\left[\langle\bm X_1^{(c)} A,\bm X_1^{(c)}\rangle\right]
   \end{equation}
   exists for each $ A\in \mathbb S(d)$.
     \end{description}
    Then  for any bounded continuous function $\varphi: D_{[0,1]}(\mathbb R^d)\to \mathbb R$,
     \begin{equation}\label{cltiid2}
    \lim_{n\to \infty} \Sbep\left[\varphi\left(\frac{\bm S_{[n\cdot]}}{\sqrt{n}}\right)\right]=\widetilde{\mathbb E}\left[\varphi(\bm W )\right],
    \end{equation}
    where $\bm W$ is a G-Brownian motion with $\bm W_1\sim N(0,G)$. In particular,
    \begin{equation}\label{cltiid3}
    \lim_{n\to \infty} \Sbep\left[\varphi\left(\frac{\bm S_n}{\sqrt{n}}\right)\right]=\widetilde{\mathbb E}\left[\varphi(\bm \xi )\right],\;\; \forall \varphi\in C_b(\mathbb R^d),
    \end{equation}
    where $\bm \xi\sim N\left(0,G\right)$.

    Conversely, if (\ref{cltiid3}) holds for any $\varphi\in C_b^1(\mathbb R^d)$ and a   random vector $\bm\xi$ with  $x^2\widetilde{\Capc}\left(|\bm \xi|\ge x\right)\to 0$ as $x\to \infty$, then
    (i-(iv) hold.

\end{theorem}

   \begin{remark} If $\Sbep[(|\bm X_1|^2-c)^+]\to 0$ as $c\to\infty$, then (i), (ii) and (iv) are satisfied, $G(A)=\Sbep\left[\langle\bm X_1 A,\bm X_1\rangle\right]$,  and (iii) is equivalent to $\Sbep[\bm X_1]=\cSbep[\bm X_1]=0$. Also, if $C_{\Capc}(|\bm X_1|^2)<\infty$, then (i), (ii) and (iv) are satisfied.

   For the one-dimensional case $d=1$, (iv) is equivalent to $\lim_{c\to \infty}\Sbep[X_1^2\wedge c]$ and  $\lim_{c\to \infty}\cSbep[X_1^2\wedge c]$ are finite which are implied by (i). In general, we don't know  whether (iv) can be derived from (i)-(iii) or not.
   \end{remark}

{\bf Proof}.   When $d=1$, this theorem is proved by Zhang\cite{Zh19b}, where it is shown that $\lim_{c\to \infty}\Sbep\left[  \bm X_1^c\right]$ and $\lim_{c\to \infty}\Sbep\left[-\bm X_1^c\right]$ exist and are finite under the condition (i).
Note
$$\left|\Sbep\left[\langle\bm X_1^c A,\bm X_1^c\rangle\right]-\Sbep\left[\langle\bm X_1^c \overline{A},\bm X_1^c\rangle\right]\right|\le |A-\overline{A}|
\Sbep[|\bm X_1|^2\wedge (d c^2)]. $$
It is easily seen that, if the limit in (\ref{cltiiid1}) exists, then it is finite and $G(A)$ is a continuous sub-linear function monotonic in $A \in  \mathbb S(d)$.  We first prove the direct part. Let $\bm Y_{n,k}= \frac{1}{\sqrt{n}}\bm X_k^{(\sqrt{n})}$. As shown in Zhang\cite{Zh19b}, by (i)-(iii) we have that
\begin{equation}\label{cltiidproof3}  \sum_{k=1}^n \Sbep[|\bm Y_{n,k}|^p]\to 0, \;\; \forall p>2,
\end{equation}
\begin{equation}\label{cltiidproof4}  \sum_{k=1}^n \left(\left|\Sbep[\bm Y_{n,k}]\right|+\left|\Sbep[-\bm Y_{nk}]\right|\right)   \to 0,
\end{equation}
\begin{equation}\label{cltiidproof5} \sum_{k=1}^n\Sbep\left[|\bm Y_{n,k}|^2\right]= \Sbep\left[|\bm X_1|^2\wedge n\right]\le C_0.
\end{equation}
  Further, by (iv),
 $$  \sum_{k=1}^{[nt]}\Sbep\left[\langle\bm Y_{n,k} A,\bm Y_{n,k}\rangle\right]=\frac{[nt]}{n} \Sbep\left[\langle\bm X_1^{(\sqrt{n})} A,\bm X_1^{(\sqrt{n})}\rangle\right]\to G(A) t. $$
 Denote $\bm W_n(t)=  \sum_{k=1}^{[nt]} \bm Y_{n,k}$. By Theorem \ref{th2}, for any bounded continuous function $\varphi: D_{[0,1]}(\mathbb R^d)\to \mathbb R$,
 \begin{equation}\label{cltiidproof7}
    \lim_{n\to \infty} \Sbep\left[\varphi\left(\bm W_n\right)\right]=\widetilde{\mathbb E}\left[\varphi(\bm W )\right].
    \end{equation}
 Note
 \begin{align}\label{cltiidproof8}\nonumber
 & \left|\Sbep\left[\varphi\left(\frac{\bm S_{[n\cdot]}}{\sqrt{n}}\right)\right]-\Sbep\left[\varphi\left(\bm W_n\right)\right]\right| \\
\le & \|\varphi\|\sum_{k=1}^n \Capc\left(\frac{\bm X_k}{\sqrt{n}}\ne \bm Y_{n,k}\right)\le \|\varphi\| n \Capc\left( |\bm X_1|\ge \sqrt{n}\right)\to 0.
\end{align}
(\ref{cltiid2}) is proved.

Now, suppose that (\ref{cltiid3}) holds. By (\ref{cltiid3}), for each element $X_{1,i}$ of $\bm X_1=(X_{1,1},\ldots, X_{1,d})$, $i=1,\ldots, d$,
we have
$$
\lim_{n\to \infty} \Sbep\left[\varphi\left(\frac{\sum_{k=1}^n X_{k,i}}{\sqrt{n}}\right)\right]=\widetilde{\mathbb E}\left[\varphi( \xi_i )\right],\;\; \forall \varphi\in C_b^1(\mathbb R),
$$
By Theorem 4.2 of Zhang\cite{Zh19b},  $\lim\limits_{c\to\infty} \Sbep[X_{1,i}^2\wedge c]$ is finite,  $x^2\Capc\left(|X_{1,i}|\ge x\right)\to 0$ as $x\to \infty$,
and  $\lim\limits_{c\to \infty}\Sbep\big[X_{1,i}^{(c)}\big]=\lim\limits_{c\to \infty}\Sbep\big[-X_{1,i}^{(c)}\big]=0$. So, (i)-(iii) are proved.

At last, we show (iv). Let $\bm Y_{n,k}$ be defined as above. Then (\ref{cltiidproof3})-(\ref{cltiidproof5}) remain true. Let $\bm T_{n,m}=\sum_{m=1}^n \bm Y_{n,m}$, $1\le m\le n$, and $\bm T_n=\bm T_{n,n}$.  Then by (\ref{eqth2proof.6}),
$$ \max_n\Sbep\Big[|\bm T_n|^p\Big]\le \max_n\Sbep\Big[\max_{m\le n}|\bm T_m|^p\Big]\le C_p, \;\; \forall p\ge 2. $$
Hence
  \begin{equation}\label{cltiidproof9}\left\{ |\bm T_n|^p; n\ge 1 \right\} \text{ is uniformly integrable for any } p\ge 2.
  \end{equation}
On the other hand, by (\ref{cltiid3}) and  (\ref{cltiidproof8}),
  \begin{equation}\label{cltiidproof9ad}
    \lim_{n\to \infty} \Sbep\left[\varphi\left(  \bm T_n\right)\right]=\widetilde{\mathbb E}\left[\varphi(\bm \xi )\right],\;\; \forall \varphi\in C_b^1(\mathbb R^d).
    \end{equation}
    Choosing $\varphi(\bm x)=|\bm x|^p\wedge c$ yields
    $$ \widetilde{\mathbb E}[|\bm \xi|^p\wedge c]=\lim_{n\to \infty} \Sbep\left[|\bm T_n |^p\wedge c\right]\le C_p. $$
  Hence
  \begin{equation}\label{cltiidproof10} \lim_{c\to\infty} \widetilde{\mathbb E}[|\bm \xi|^p\wedge c] \le C_p\text{ is finite for any } p\ge 2.
  \end{equation}
  Let $G_{\xi}(A,c)=\widetilde{\mathbb E}\left[\langle\bm \xi^{(c)}A,\bm \xi^{(c)}\rangle\right]$. Note, for $a>b$,
  \begin{align*}
  \left|\langle\bm \xi^{(a)}A,\bm \xi^{(b)}\rangle- \langle\bm \xi^{(a)}A,\bm \xi^{(b)}\rangle\right|
  \le |A|(|\bm \xi^{(a)}|+|\bm \xi^{(b)}|)|\bm \xi^{(a)}-\bm \xi^{(b)}|.
  \end{align*}
  It follows that
  \begin{align}\label{cltiidproof11} &|G_{\xi}(A,a)-G_{\xi}(A,b)|\nonumber \\
  \le & |A|\Big(\widetilde{\mathbb E}\big[(|\bm \xi^{(a)}|+|\bm \xi^{(b)}|)\big]\Big)^{1/2}
  \Big(\widetilde{\mathbb E}\big[\sum_{k=1}^d (\xi_k^2\wedge a^2-b^2)^+\big]\Big)^{1/2}\\
  \le & C|A|\left(\Sbep\big[|\bm \xi|^2\wedge (d a^2)\big]\right)^{1/2} \left(\frac{d \Sbep[|\bm \xi|^3\wedge a^3]}{b}\right)^{1/2} \nonumber\\
  \to & 0\text{ as } a>b\to \infty,\nonumber
   \end{align}
   by (\ref{cltiidproof10}). If follows that
   $$ G_{\xi}(A)=\lim_{c\to \infty} G_{\xi}(A,c) \text{ exists and is finite}. $$
   Now, choosing $\varphi(\bm x)=\langle \bm x^{(c)} A,\bm x^{(c)}\rangle$ in (\ref{cltiidproof9ad}) yields
   $$ \lim_{n\to \infty} \Sbep\Big[\langle \bm T_n^{(c)} A,\bm T_n^{(c)}\rangle\Big]=G_{\xi}(A, c). $$
   Note that $|\langle \bm T_n A,\bm T_n- \langle \bm T_n^{(c)} A,\bm T_n^{(c)}\rangle|\le 2|A|\cdot |\bm T_n|^2I\{|\bm T_n|>c\}$, and $\{|\bm T_n|^2, n\ge 1\}$ is uniformly integrable by (\ref{cltiidproof9}). Letting $c\to \infty$ in the above equation yields
   $$ \lim_{n\to \infty} \Sbep\big[\langle \bm T_n A,\bm T_n\rangle\big]=G_{\xi}(A). $$
   On the other hand, note
   \begin{align*}
   \langle \bm T_n A,\bm T_n\rangle=\sum_{k=1}^n\langle \bm Y_{n,k} A,\bm Y_{n,k}\rangle+2\sum_{k=1}^n\langle \bm T_{n,k-1} A,\bm Y_{n,k}\rangle.
   \end{align*}
   Since
   $$\Sbep\left[\langle \bm x,\bm X\rangle\right]\le \sum_{i=1}^d( x_i^+ \Sbep[X_i]+x_i^-\Sbep[-X_i])\le 2 |\bm x| \big(|\Sbep[\bm Y_{n,k}]|+|\Sbep[-\bm Y_{n,k}]|\big), $$
   we have
    \begin{align*}
  \Sbep\left[\pm \sum_{k=1}^n\langle \bm T_{n,k-1} A,\bm Y_{n,k}\rangle\right]
  \le  &2 \sum_{k=1}^n \Sbep[|\bm T_{n,k-1} A|]\big(|\Sbep[\bm Y_{n,k}]|+|\Sbep[-\bm Y_{n,k}]|\big)\\
  \le & C\sum_{k=1}^n \big(|\Sbep[\bm Y_{n,k}]|+|\Sbep[-\bm Y_{n,k}]|\big)\to 0.
   \end{align*}
  It follows that
    \begin{align*} &\Sbep\big[\langle \bm T_n A,\bm T_n\rangle\big]-\sum_{k=1}^n \Sbep\big[\langle \bm Y_{n,k} A,\bm Y_{n,k}\rangle\big]\\
    = &\Sbep\big[\langle \bm T_n A,\bm T_n\rangle\big]-  \Sbep\big[\sum_{k=1}^n\langle \bm Y_{n,k} A,\bm Y_{n,k}\rangle\big]\to 0.
    \end{align*}
   We conclude that
   $$ \Sbep\big[\langle \bm X_1^{(\sqrt{n})} A,\bm X_1^{(\sqrt{n})}\rangle\big]=\sum_{k=1}^n \Sbep\big[\langle \bm Y_{n,k} A,\bm Y_{n,k}\rangle\big]
   \to G_{\xi}(A). $$
   Similar to (\ref{cltiidproof11}), for $\sqrt{n}\le b\le a\le \sqrt{n+1}$ we have
   \begin{align*}  &|\Sbep\big[\langle \bm X_1^{(a)} A,\bm X_1^{(a)}\rangle\big]-\Sbep\big[\langle \bm X_1^{(b)} A,\bm X_1^{(b)}\rangle\big]| \\
  \le & |A|\Big(\widetilde{\mathbb E}\big[(|\bm X_1^{(a)}|+|\bm X_1^{(b)}|)\big]\Big)^{1/2}
  \Big(\widetilde{\mathbb E}\big[\sum_{k=1}^d (X_{1,k}^2\wedge a^2-b^2)^+\big]\Big)^{1/2}\\
  \le & C |A| \Big(\sum_{k=1}^n (n+1)\widetilde{\Capc}\big(|X_{1,k}|\ge \sqrt{n}\big)\Big)^{1/2}\to 0,
   \end{align*}
   by (i) and (iii). Hence,
   $$\lim_{c\to \infty}\Sbep\big[\langle \bm X_1^{(c)} A,\bm X_1^{(c)}\rangle\big]=\lim_{n\to\infty}
   \Sbep\big[\langle \bm X_1^{(\sqrt{n})} A,\bm X_1^{(\sqrt{n})}\rangle\big]=G_{\xi}(A), \;\; A\in \mathbb S(d). $$
   (iv) is now proved. \hfill $\Box$.

\bigskip

At last,  we   give a L\'evy   characterization of a multi-dimensional G-Brownian motion as an application of Theorem \ref{thFCLTM}.    Let $\{\mathscr{H}_t; t\ge 0\}$   be  a non-decreasing family of subspaces of $\mathscr{H}$ such that (1) a constant $c\in \mathscr{H}_t$ and, (2) if $X_1,\ldots,X_d\in \mathscr{H}_t$, then $\varphi(X_1,\ldots,X_d)\in \mathscr{H}_t$ for any $\varphi\in C_{l,lip}$.
 We consider a system of operators in $\mathscr{L}(\mathscr{H})=\{X\in \mathscr{H}; \Sbep[|X|]<\infty\}$,
 $$\Sbep_t: \mathscr{L}(\mathscr{H})\to \mathscr{L}(\mathscr{H}_t) $$
 and denote $\Sbep[X|\mathscr{H}_t]=\Sbep_t[X]$, $\cSbep[X|\mathscr{H}_t]=-\Sbep_t[-X]$.
Suppose that the operators $\Sbep_t$ satisfy  the following properties:  for all $X, Y \in \mathscr{L}({\mathscr H})$,
\begin{description}
  \item[\rm (i)]   $ \Sbep_t [ X+Y]=X+\Sbep_t[Y]$ in $L_1$ if $X\in \mathscr{H}_t$, and $ \Sbep_t [ XY]=X^+\Sbep_t[Y]+X^-\Sbep_t[-Y]$ in $L_1$
if $X\in \mathscr{H}_t$ and $XY\in \mathscr{L}({\mathscr H})$;
\item[\rm (ii)]   $\Sbep\left[\Sbep_t [ X]\right]=\Sbep[X]$.
 \end{description}
For a random vector $\bm X=(X_1,\ldots, X_d)$, we denote $\Sbep_t[\bm X]=\big(\Sbep_t[X_1],\ldots, \Sbep_t[X_d]\big)$.

 \begin{definition} A $d$-dimensional process $\bm M_t$ is called a martingale, if $\bm M_t\in \mathscr{L}(\mathscr{H}_t)$  and
 $$ \Sbep[\bm M_t|\mathscr{H}_s]=\bm M_s, \;\; s\le t. $$
 \end{definition}
 Denote
 \begin{align*}
 W_T(\bm M,\delta) = &\sup_{t_i}\Sbep\big[\max_{1\le i\le n} |\bm M(t_i)-\bm M(t_{i-1})|\wedge 1\big], \\
 & \text{ where the supermum }  \sup_{t_i} \text{ is taken over all } t_is \text{ with } \\
 & 0=t_0<t_1<\ldots<t_n=T,\;\; \delta/2<t_i-t_{i-1}<\delta, \; i=1,\ldots, n.
 \end{align*}

 The L\'evy characterization of a one-dimensional $G$-Brownian motion under G-expectation   in a Wiener space is   established by Xu and Zhang \cite{XuZhang09,XuZhang10} and extended by Lin\cite{Lin13} by the method of  the stochastic calculus. The following theorem gives a L\'evy characterization of a $d$-dimensional  G-Brownian motion.
  \begin{theorem} \label{th4.2} Let $\bm M_t$ be a $d$-dimensional random process in $(\Omega,\mathscr{H},\mathscr{H}_t, \Sbep)$ with $\bm M_0=\bm 0$,
 \begin{equation} \label{eqth4.2.1} \text{ for all }  p>0  \text{ and } t\ge 0, \;\; C_{\Capc}(|\bm M_t|^p)<\infty \implies \Sbep[|\bm M_t|^p]<\infty.
 \end{equation}    Suppose that $\bm M_t$ satisfies
 \begin{description}
 \item[\rm (I)] Both $\bm M_t$ and $-\bm M_t$ are martingales;
 \item[\rm (II)] There is a  a continuous sub-linear function $G:\mathbb S(d) \to \mathbb R$     monotonic in $A \in  \mathbb S(d)$ such that $\langle \bm M_t A, \bm M_t\rangle-G(A)t$ is a real martingale
 for each $A \in  \mathbb S(d)$;
  \item[\rm (III)]  For any $T>0$, $\lim_{\delta\to 0}W_T(\bm M,\delta)=0$.
 \end{description}
 Then, $\bm M_t$ satisfies Property (ii) as in Definition \ref{DefG-B}    with $\bm M_1\sim N(0,G)$.
\end{theorem}
{\bf Proof.} The proof is very similar to that of Theorem 5.3 of Zhang\cite{Zh19} by applying Theorem (\ref{thFCLTM}) and so is omitted. \hfill $\Box$

  \section{Proofs.}\label{sectProof}
  \setcounter{equation}{0}

To prove functional central limit theorems, we need the following  Rosenthal-type inequalities which can be proved by the same argument as in Theorem 4.1 of Zhang\cite{Zh19}.

\begin{lemma} \label{lemRosenIeq}   Suppose that $\{X_{n,i}\}$ are a set of bounded random variables, $X_{n,k}\in \mathscr{H}_{n,k}$.  Set $S_0=0$, $S_k=\sum_{i=1}^k  X_{n,i}$.  Then,
\begin{align}\label{eqlemRosenIeq.1}
\Sbep\Big[\Big(\max_{k\le {k_n}} (S_{k_n}-S_k)\Big)^2\Big|\mathscr{H}_{n,0}\Big]\le      \Sbep\Big[\sum_{k=1}^{k_n}   \Sbep[X_{n,k}^2|\mathscr{H}_{n,k-1}]\Big|\mathscr{H}_{n,0}\Big] \text{ in } L_1,
\end{align}
when $\Sbep[X_{n,k}|\mathscr{H}_{n,k-1}]\le 0$, $k=1,\ldots, k_n$. In general,   for $p\ge 2$ there is a constant $C_p$ such that
\begin{align} \label{eqlemRosenIeq.4}
& \Sbep\Big[\max_{k\le {k_n}} |S_k|^p\big|\mathscr{H}_{n,0}\Big] \nonumber \\
&\le    C_p\left\{\Sbep\left[\sum_{k=1}^{k_n} \Sbep[|X_{n,k}|^p|\mathscr{H}_{n,k-1}]\Big|\mathscr{H}_{n,0}\right]+\Sbep\left[\Big(\sum_{k=1}^{k_n}   \Sbep[X_{n,k}^2|\mathscr{H}_{n,k}]\Big)^{p/2}\Big|\mathscr{H}_{n,0}\right]\right.\nonumber\\
&  \qquad \left.+\Sbep\left[\Big\{\sum_{k=1}^{k_n}  \Big(\big( \Sbep[X_{n,k} |\mathscr{H}_{n,k}]\big)^++\big( \cSbep[X_{n,k} |\mathscr{H}_{n,k}]\big)^-\Big)\Big\}^p\Big|\mathscr{H}_{n,0}\right]\right\} \text{ in } L_1.
\end{align}
\end{lemma}

 {\bf Proof of Theorem \ref{thFCLTM}}. With the same arguments as those  in  the proofs of Theorems 3.1 and 3.2 of Zhang\cite{Zh19},   we can assume that $\delta_{k_n}=\sum_{k=1}^{k_n}\Sbep[|\bm Z_{n,k}|^2|\mathscr{H}_{n,k-1}]\le 2\rho(1) $ in $L_1$, $\chi_{k_n}=:\sum_{k=1}^{k_n}\left\{ |\Sbep[ \bm Z_{n,k} |\mathscr{H}_{n,k-1}]|+|\cSbep[\bm Z_{n,k} |\mathscr{H}_{n,k-1}]|\right\}<1$ in $L_1$ and $|\bm Z_{n,k}|\le \epsilon_n$, $k=1,\ldots,k_n$, with a sequence $0<\epsilon_n\to 0$.  Under these assumptions,  the property (g) of the conditional expectation implies  that  all random variables considered above  are bounded in $L_p$ for all $p>0$, and then the convergences in (\ref{eqCLTCondM.2}) and (\ref{eqFCLTCondM.3}) all hold in $L_p$ for any $p>0$, by Lemma \ref{lemma4.0.1}.   As in Zhang\cite{Zh19},   it can be  shown that for any $\epsilon>0$,
\begin{equation}\label{eqpoofthFCLTMtight} \lim_{\delta\to 0} \limsup_{n\to \infty} \Capc\left( w_{\delta}\left(\bm W_n\right)\ge  \epsilon\right)=0,
\end{equation}
where $\omega_{\delta}(\bm x)=\sup_{|t-s|<\delta,t,s\in[0,1]}|\bm x(t)-\bm x(s)|$. So, for (\ref{eqFCLTM}) it is sufficient to show (\ref{eqfinitedimension}). With the same argument of Zhang\cite{Zh19}, it is sufficient to show that for any $0\le s<t\le 1$ and  a bounded Lipschitz   function $\varphi(\bm u, \bm x)$,
 \begin{equation}\label{eqproofFCLTM.2}
\Sbep\left[\left|\Sbep\left[\varphi\big(\bm u, \bm S_{n, \tau_n(t)}-\bm S_{n, \tau_n(s)}\big)\big|\mathscr{H}_{n,\tau_n(s)}\right]
-\widetilde{\mathbb E}\left[\varphi\big(\bm u,  \bm W(\rho(t))-\bm W(\rho(s)) \big)\right]\right|\right]\to 0.
\end{equation}

We first show that, for any $r\ge 2$ there is a positive constant $C_r>0$ such that
\begin{align}
\label{eq5.17}
 &\Sbep\left[\max_{\tau_n(s)\le k \le \tau_n(t)}  \left|\bm S_{n,k}-\bm S_{n,\tau_n(s)}\right|^r\big|\mathscr{H}_{n,\tau_n(s)}\right]\le  C_r \; \text{ in } L_p, \\
 \label{eq5.18}
&\Sbep\left[\left|\bm S_{n,\tau_n(t)}-\bm S_{n,\tau_n(s)}\right|^r\big|\mathscr{H}_{n,\tau_n(s)}\right]\le  C_r\left(\rho(t)-\rho(s)\right)^{p/2}+o(1) \; \text{ in } L_p,  \\
\label{eq5.19}
& \Sbep\left[\bm S_{n,\tau_n(t)}-\bm S_{n,\tau_n(s)} \big|\mathscr{H}_{n,\tau_n(s)}\right]\to  \bm 0   \; \text{ in } L_p, \\
 \label{eq5.20}
&  \cSbep\left[\bm S_{n,\tau_n(t)}-\bm S_{n,\tau_n(s)} \big|\mathscr{H}_{n,\tau_n(s)}\right]\to  \bm 0 \; \text{ in } L_p,
\\
\label{eq5.21}
& \Sbep\left[ \Big\langle(\bm S_{n,\tau_n(t)}-\bm S_{n,\tau_n(s)})A, \bm S_{n,\tau_n(t)}-\bm S_{n,\tau_n(s)}\Big\rangle\Big|\mathscr{H}_{n,\tau_n(s)}\right]\nonumber \\
& \quad \to G(A) \big(\rho(t)-\rho(s)\big) \;  \text{ in } L_p, \;\; \forall A \in \mathbb S(d),
\end{align}
for any $0<s<t$ and $p>0$. Further, (\ref{eq5.21}) holds uniformly in $A \in \mathbb S(d)$ with $|A|\le c$.

For (\ref{eq5.17})-(\ref{eq5.20}), it is sufficient to verify the one-dimensional case. For (\ref{eq5.17}), by Lemma \ref{lemRosenIeq},
\begin{align}
&\Sbep\left[\max_{\tau_n(s)\le k \le\tau_n(t)} \left|S_{n,k}-S_{n,\tau_n(s)}\right|^r\big|\mathscr{H}_{n,\tau_n(s)}\right] \nonumber\\
\le & C_r  \left\{\Sbep\left[\sum_{k=\tau_n(s)+1}^{\tau_n(t)} \Sbep[|Z_{n,k}|^r|\mathscr{H}_{n,k-1}]\Big|\mathscr{H}_{n,\tau_n(s)}\right] \right. \nonumber\\
& + \Sbep\left[\Big(\sum_{k=\tau_n(s)+1}^{\tau_n(t)}  \Sbep[|Z_{n,k}|^2|\mathscr{H}_{n,k}]\Big)^{r/2}\Big|\mathscr{H}_{n,\tau_n(s)}\right] \nonumber\\
&    +\left. \Sbep\left[\Big\{\sum_{k=\tau_n(s)+1}^{\tau_n(t)}  \Big(\big| \Sbep[Z_{n,k} |\mathscr{H}_{n,k}]\big|+\big| \cSbep[Z_{n,k} |\mathscr{H}_{n,k}]\big|\Big)\Big\}^r\Big|\mathscr{H}_{n,\tau_n(s)}\right]\right\}\nonumber  \\
\le & C_r  \left\{\epsilon_n^{r-2} \Sbep\left[\delta_{k_n} \Big|\mathscr{H}_{n,\tau_n(s)}\right] + \Sbep\left[\Big(\sum_{k=\tau_n(s)+1}^{\tau_n(t)}  \Sbep[ |Z_{n,k}|^2|\mathscr{H}_{n,k}]\Big)^{r/2}\Big|\mathscr{H}_{n,\tau_n(s)}\right]\right. \nonumber\\
&    +\left. \Sbep\left[\chi_{k_n}^r\Big|\mathscr{H}_{n,\tau_n(s)}\right]\right\}\label{eq5.23}  \\
\le & C_r\left\{ 2\rho(1)+(2\rho(1))^{r/2}+1\right\} \text{ in } L_1. \nonumber
\end{align}
Note that the random variable $\max\limits_{\tau_n(s)\le k \le\tau_n(t)}  |S_{n,k}-S_{n,\tau_n(s)}|$ is a bounded ($\le (\tau_n(t)-\tau_n(s))\epsilon_n$). By the property (g) of $\Sbep_{n,k}$, $\Sbep\Big[\max\limits_{\tau_n(s)\le k \le\tau_n(t)}  |S_{n,k}-S_{n,\tau_n(s)}|^r\big|\mathscr{H}_{n,\tau_n(s)}\Big]$ is bounded in $L_p$ for any $p>0$. Hence, by (1) and (2) of Lemma \ref{lemma4.0.1}, (\ref{eq5.17}) is proved.  By this inequality and Lemma \ref{lemma4.0.1}, it is sufficient to consider the case of $p=1$ for (\ref{eq5.18})-(\ref{eq5.21}).

It is easily shown that
\begin{align*}
\Sbep\left[\pm\left(S_{n,\tau_n(t)}-S_{n,\tau_n(s)}\right) \big|\mathscr{H}_{n,\tau_n(s)}\right]\le \Sbep\left[\chi_{k_n}\big|\mathscr{H}_{n,\tau_n(s)}\right]\to 0 \text{ in } L_1,
\end{align*}
which implies (\ref{eq5.19}) and (\ref{eq5.20}).

For (\ref{eq5.21}), we first note that
\begin{equation}
 \label{eq5.24} \sum_{k=\tau_n(s)+1}^{\tau_n(t)}\Sbep\left[\big\langle \bm Z_{n, k}A,\bm Z_{n, k}\big\rangle\Big|\mathscr{H}_{n,k-1}\right]\to G(A)\big(\rho (t)-\rho(s)\big) \text{ in } L_p, \\
 \end{equation}
for any $p>0$, by condition (\ref{eqFCLTCondM.3}).  Without loss of generality, we assume $s=0$, $t=1$. Note
\begin{align*}
&\langle \bm S_{k_n}A, \bm S_{k_n}\rangle -\sum_{k=1}^{k_n}\Sbep\left[\langle \bm Z_{n, k}A, \bm Z_{n,k}\rangle\big|\mathscr{H}_{n,k-1}\right] \\
=& \sum_{k=1}^{k_n}\left(\langle \bm Z_{n, k}A, \bm Z_{n,k}\rangle-\Sbep\left[\langle \bm Z_{n, k}A, \bm Z_{n,k}\rangle\Big|\mathscr{H}_{n,k-1}\right]\right)+2\sum_{k=1}^{k_n}\langle \bm S_{n, k-1}A, \bm Z_{n,k}\rangle,
\end{align*}
\begin{align*}
&\Sbep\left[\pm \langle \bm S_{n, k-1}A, \bm Z_{n,k}\rangle\big|\mathscr{H}_{n,k-1}\right]\\
\le & 2|\bm S_{n,k-1}A| \left\{\left|\Sbep\left[ \bm Z_{n,k}\big|\mathscr{H}_{n,k-1}\right]\right|+\left|\Sbep\left[- \bm Z_{n,k}\big|\mathscr{H}_{n,k-1}\right]\right|\right\}
\text{ in } L_1.
\end{align*}
And then
\begin{align*}
&\Sbep\left[\pm \left(\sum_{k=1}^{k_n} \langle \bm S_{n, k-1}A, \bm Z_{n,k}\rangle\right)\big|\mathscr{H}_{n,0}\right] \\
\le & 2\Sbep\left[ \sum_{k=1}^{k_n}|\bm S_{n,k-1}A|\left\{ \left|\Sbep\left[ \bm Z_{n,k}\big|\mathscr{H}_{n,k-1}\right]\right|
+\left|\cSbep\left[ \bm Z_{n,k}\big|\mathscr{H}_{n,k-1}\right]\right|\right\}\big|\mathscr{H}_{n,0}\right]\text{ in } L_1.
\end{align*}
It follows that
\begin{align*}
& \left|\Sbep\Big[\langle \bm S_{k_n}A, \bm S_{k_n}\rangle -\sum_{k=1}^{k_n}\Sbep\left[\langle \bm Z_{n, k}A, \bm Z_{n,k}\rangle\big|\mathscr{H}_{n,k-1}\right] \Big|\mathscr{H}_{n,0}\Big]\right|
\le 2\Sbep\left[  \chi_{k_n}\max_{k\le k_n}|\bm S_{n,k}A|    \Big|\mathscr{H}_{n,0}\right]\text{ in } L_1.
\end{align*}
Taking the sub-linear expectation yields
\begin{align*}
&\Sbep\left[\left|\Sbep\Big[\langle \bm S_{k_n}A, \bm S_{k_n}\rangle -\sum_{k=1}^{k_n}\Sbep\left[\langle \bm Z_{n, k}A, \bm Z_{n,k}\rangle\big|\mathscr{H}_{n,k-1}\right] \Big|\mathscr{H}_{n,0}\Big]\right|\right] \\
\le &2\Sbep\Big[  \chi_{k_n}\max_{k\le k_n}|\bm S_{n,k}A|\Big] \le 2\left(\Sbep[\chi_{k_n}^2]\Sbep[\max_{k\le k_n}|\bm S_{n,k}A|^2]\right)^{1/2}
\le   C   \left(\Sbep[\chi_{k_n}^2]\right)^{1/2}\to 0,
\end{align*}
by (\ref{eq5.17}) and the fact that $\chi_{k_n}\to 0$ in $L_p$. By noting (\ref{eq5.24}), we have
$$ \Sbep\left[\left|\Sbep\Big[\langle \bm S_{k_n}A, \bm S_{k_n}\rangle -G(A)\rho(1) \Big|\mathscr{H}_{n,0}\Big]\right|\right]\to 0. $$
(\ref{eq5.21}) is proved.  By the same argument as in Remark \ref{remark2.1},  (\ref{eq5.21}) holds uniformly in $A\in \mathbb S(d)$ with $|A|\le c$.

For (\ref{eq5.18}), it is easily seen the first and the third terms in (\ref{eq5.23}) converges to 0 in $L_1$, and the second term converges to $\big(\rho(t)-\rho(s)\big)^{r/2}$  by (\ref{eq5.24}). And hence, (\ref{eq5.18}) is proved.

Now, we tend to prove (\ref{eqproofFCLTM.2}). Without loss of generality, we assume $s=0$ and $t=1$.
 Let $V(t, \bm x)=V^{\bm u} (t, \bm x)$ be the unique viscosity solution of the following equation,
$$ \partial_t V^{\bm u} + \frac{1}{2}G( D^2 V^{\bm u})=0,\;\;  (t, x) \in [0,\varrho+ h] \times \mathbb R, \; V^{\bm u}|_{t=\varrho+h} = \varphi(\bm u, \bm x),
$$
where $\varrho=\rho(1)-\rho(0)$. Without loss of generality, we assume that there is a constant $\epsilon>0$ such that
\begin{equation}\label{eqInterRegul.1}
 G(A)-G(\overline{A})\ge tr(A-\overline{A}) \epsilon \text{ for all } \; A,\overline{A}\in \mathbb S(d) \text{ with } A\ge \overline{A},
 \end{equation} for otherwise we can add a random vector $\epsilon\Sbep[|\bm Z_{n,k}|^2\big|\mathscr{H}_{n,k-1}]\bm \xi_{n,k}$ to $\bm Z_{n,k}$, where $\bm \xi_{n,k}$  has a $d$-dimensional standard normal $N(0, I_{d\times d})$ distribution and is independent to $\bm Z_{n,1},\ldots, \bm Z_{n,k}$, $\bm \xi_{n,1},\ldots,\bm \xi_{n,k-1}$. Under (\ref{eqInterRegul.1}) by
the interior regularity of $V^{\bm u}$ (c.f. Theorem 4.5 of Peng\cite{Peng10}),
\begin{equation}\label{eqPDE2}\|V^{\bm u}\|_{C^{1+\alpha/2,2+\alpha}([0,\rho+h/2]\times \mathbb R^d)} < \infty, \text{ for some } \alpha\in (0, 1).
\end{equation}
According to the definition of $G$-normal distribution, we have $V^{\bm u}(t,\bm x)=\widetilde{\mathbb E}\big[\varphi(\bm u, \bm x+\sqrt{\varrho+h-t}\bm \xi)\big]$ where $\bm\xi\sim N(0,G)$ under $\widetilde{\mathbb E}$. In particular,
$$ V^{\bm u} (h,\bm 0)=\widetilde{\mathbb E}\big[\varphi( \bm u, \sqrt{\varrho}\bm\xi)\big], \;\; V^{\bm u}(\varrho+h, \bm x)=\varphi(\bm u,\bm x). $$
 Following the proof of Theorem 3.1 and 3.2 of Zhang\cite{Zh19}, it is sufficient to show that
\begin{equation}\label{eqconV}
 \Sbep\left[\Big|\Sbep[ V(\varrho, \bm S_{k_n}) |\mathscr{H}_{n,0}]-V(0,\bm 0)\Big|\right]\to 0.
 \end{equation}
As in Zhang\cite{Zh19}, it can be  proved that,  for all $(t,\bm x)\in [0, \varrho+h/2] \times \mathbb R^d$,
$$ |D V(t,\bm x)|\le C,\;\;  |\partial_t V(t,\bm x)|\le C, \;\;  |D^2 V(t, \bm x)|\le   C+C|\bm x|^{\alpha}. $$
For an integer $m$ large enough, we define $t_i=i/m$, $\bm Y_{n,i}=\bm S_{n,\tau_n(t_i)}-\bm S_{n,\tau_n(t_{i-1})}$, $\widetilde{\delta}_i=\rho(t_i)$, $\bm T_i=\sum_{j=1}^i \bm Y_{n,j}$, $i=1,\ldots, m$.
Applying the Taylor's expansion yields
\begin{align*}
 & V(\varrho, \bm S_{k_n})-V(0,\bm 0)\\
 =&\sum_{i=0}^{m-1} \left\{[ V( \widetilde{\delta}_{i+1},  \bm T_{i+1})-V( \widetilde{\delta}_i, \bm T_{i+1})]
  +[ V( \widetilde{\delta}_i,  \bm T_{i+1})-V( \widetilde{\delta}_i,  \bm T_i)]\right\} \\
 =: & \sum_{i=0}^{m-1} \left\{I_{n}^i+J_{n}^i\right\},
\end{align*}
with
\begin{align*}
 J_{n}^i =&  \partial_t V( \widetilde{\delta}_i,  \bm T_i)\big(\widetilde{\delta}_{i+1}-\widetilde{\delta}_i\big)+
\big\langle D V( \widetilde{\delta}_i,  \bm T_i), \bm Y_{n,i+1}\big\rangle     +\frac{1}{2}\big\langle \bm Y_{n,i+1} D^2 V(\widetilde{\delta}_i,  \bm T_i), \bm Y_{n,i+1}\big\rangle          \\
=&   \left\{\partial_t V( \widetilde{\delta}_i,  \bm T_i)+\frac{1}{2}G\Big(D^2 V(\widetilde{\delta}_i,  \bm T_i)\Big)\right\}\big(\widetilde{\delta}_{i+1}-\widetilde{\delta}_i\big) \\
  &+ \frac{1}{2} \Big\{ \big\langle \bm Y_{n,i+1} D^2 V(\widetilde{\delta}_i,  \bm T_i), \bm Y_{n,i+1}\big\rangle-\Sbep\left[\big\langle \bm Y_{n,i+1} D^2 V(\widetilde{\delta}_i,  \bm T_i), \bm Y_{n,i+1}\big\rangle\Big|\mathscr{H}_{n,\tau_n(t_i)}\right]\Big\}  \\
&  + \Big\{\big\langle D V( \widetilde{\delta}_i,  \bm T_i), \bm Y_{n,i+1}\big\rangle \Big\}\\
 &  +\frac{1}{2} \Big\{ \Sbep\left[\big\langle \bm Y_{n,i+1} D^2 V(\widetilde{\delta}_i,  \bm T_i), \bm Y_{n,i+1}\big\rangle\Big|\mathscr{H}_{n,\tau_n(t_i)}\right]
   -G\Big(D^2 V(\widetilde{\delta}_i,  \bm T_i)\Big)\big(\widetilde{\delta}_{i+1}-\widetilde{\delta}_i\big)\Big\}
   \\
 =:&0+ J_{n,1}^i+J_{n,2}^i+J_{n,3}^i
\end{align*}
 and
\begin{align*}
I_{n}^i=&\big(\widetilde{\delta}_{i+1}-\widetilde{\delta}_i\big)\left[  \big( \partial_t V( \widetilde{\delta}_i+\gamma \big(\widetilde{\delta}_{i+1}-\widetilde{\delta}_i\big),  \bm T_{i+1})-\partial_t V( \widetilde{\delta}_i ,  \bm T_{i+1})\big) \right.\\
&\left.\qquad \quad +\big( \partial_t V( \widetilde{\delta}_i,  \bm T_{i+1})-\partial_t V( \widetilde{\delta}_i ,  \bm T_i)\big)   \right] \\
&+ \frac{1}{2}\left\langle\bm Y_{n,i+1}\left[ D^2 V(\widetilde{\delta}_i,  \bm T_i+\beta \bm Y_{n,i+1})-D^2 V(\widetilde{\delta}_i,  \bm T_i)\right],\bm Y_{n,i+1}\right\rangle,
\end{align*}
where $\gamma$ and $\beta$ are between $0$ and $1$.

By (\ref{eqPDE2}), it is easily seen that
\begin{align*}
|I_n^i|\le  &  C\big| \widetilde{\delta}_{i+1}-\widetilde{\delta}_i \big|^{2+\alpha } +C (\widetilde{\delta}_{i+1}-\widetilde{\delta}_i) |\bm Y_{n,i+1}|^{\alpha} +C|\bm Y_{n,i+1}|^{2+\alpha}\\
\le & C \left(\rho(t_{i+1})-\rho(t_i)\right)^{1+\alpha/2} +o(1)\; \text{ in } L_1,
\end{align*}
by (\ref{eq5.18}), where $C$ is a positive constant which does not depend on $t_i$s.

For $J_{n,1}^i$,   it follows that
$$ \Sbep\left[J_{n,1}^i\big|\mathscr{H}_{n,\tau_n(t_i)}\right] =0\;\text{ in } L_1. $$
It follows that
\begin{align*}\label{eqtheorem2proof.3}
\Sbep\left[\sum_{i=0}^{m-1}J_{n,1}^i|\mathscr{H}_{n,0}\right]
=& \Sbep\left[\sum_{i=0}^{m-2}J_{n,1}^i+\Sbep\left[J_{n,1}^{m-1}\big|\mathscr{H}_{n,\tau_n(t_{m-1})}\right]\Big|\mathscr{H}_{n,0}\right]\\
=
&\Sbep\left[\sum_{i=0}^{m-2}J_{n,1}^i\Big|\mathscr{H}_{n,0}\right]=\ldots=0 \text{ in } L_1.
\end{align*}

For $J_{n,2}^i$, we have
\begin{align*}
&\Sbep[J_{n,2}^i|\mathscr{H}_{n,0}]=\Sbep\Big[\Sbep[J_{n,2}^i|\mathscr{H}_{n,\tau_n(t_i)}]\big|\mathscr{H}_{n,0}\Big]\\
\le & \Sbep\left[ |D V( \widetilde{\delta}_i,  T_i))|\left\{ | \Sbep[\bm Y_{n,i+1}|\mathscr{H}_{n,\tau_n(t_i)}]|+ |\cSbep[\bm Y_{n,i+1}|\mathscr{H}_{n,\tau_n(t_i)}]|\right\} \Big| \mathscr{H}_{n,0}\right] \\
\le &C \Sbep\left[  \left|\Sbep[\bm Y_{n,i+1}|\mathscr{H}_{n,\tau_n(t_i)}]\right|+ \left|\cSbep[\bm Y_{n,i+1}|\mathscr{H}_{n,\tau_n(t_i)}]\right|  \Big| \mathscr{H}_{n,0}\right]\to 0 \text{ in } L_1,
 \end{align*}
by (\ref{eq5.19}) and (\ref{eq5.20}). Similarly,
$\Sbep[-J_{n,2}^i|\mathscr{H}_{n,0}]\le o(1)$   in $L_1$.

For $J_{n,3}^i$, we have
\begin{align*}
|J_{n,3}^i|\le & \frac{1}{2} |D^2 V(\widetilde{\delta}_i,  \bm T_i)|\sup_{|A|\le 1}\Big| \Sbep\left[\big\langle \bm Y_{n,i+1} A, \bm Y_{n,i+1}\big\rangle\Big|\mathscr{H}_{n,\tau_n(t_i)}\right]
   -G(A)\big(\widetilde{\delta}_{i+1}-\widetilde{\delta}_i\big)\Big|\\
   \le & (C+C|\bm T_i|^{\alpha})\sup_{|A|\le 1}\Big| \Sbep\left[\big\langle \bm Y_{n,i+1} A, \bm Y_{n,i+1}\big\rangle\Big|\mathscr{H}_{n,\tau_n(t_i)}\right]
   -G(A)\big(\widetilde{\delta}_{i+1}-\widetilde{\delta}_i\big)\Big| \\
\le &C \sup_{|A|\le 1}\Big| \Sbep\left[\big\langle \bm Y_{n,i+1} A, \bm Y_{n,i+1}\big\rangle\Big|\mathscr{H}_{n,\tau_n(t_i)}\right]
   -G(A)\big(\widetilde{\delta}_{i+1}-\widetilde{\delta}_i\big)\Big|=o(1) \text{ in } L_1
\end{align*}
by (\ref{eq5.17}) and (\ref{eq5.21}), where $C$ is a positive constant which does not depend on $t_i$s.

Combing the above arguments yields
\begin{align*}
 &\Big|\Sbep[ V(\varrho, \bm S_{k_n}) |\mathscr{H}_{n,0}]-V(0,\bm 0)\Big| \\
 \le & \sum_{i=0}^{m-1} \left\{ \Sbep[J_{n,2}^i|\mathscr{H}_{n,0}]+\Sbep[-J_{n,2}^i|\mathscr{H}_{n,0}]\right. \\
 & \qquad \left. +\Sbep[|J_{n,3}^i|\big|\mathscr{H}_{n,0}]
 +\Sbep[|I_n^i|\big|\mathscr{H}_{n,0}]\right\}\\
 \le & C \sum_{i=0}^{m-1}\left(\rho(t_{i+1})-\rho(t_i)\right)^{1+\alpha/2} +o(1) \\
  \le & C  \max_i\Big(\rho((i+1)/m)-\rho(i/m)\Big)^{\alpha/2}\varrho +o(1)\; \text{ in } L_1.
 \end{align*}
The proof of (\ref{eqconV}) is completed by letting $m\to \infty$.  \hfill $\Box$



 \end{document}